\documentclass[12pt,reqno]{amsart}
\usepackage{amssymb}
\usepackage{amscd}
\usepackage{hyperref}

\newtheorem{theorem}{Theorem}[section]

\newtheorem{lemma}[theorem]{Lemma}
\newtheorem{proposition}[theorem]{Proposition}

\newtheorem{remark}[theorem]{Remark}

\newtheorem{example}[theorem]{Example}
\numberwithin{equation}{section}
\begin{document}
\title{Natural equivariant Dirac operators}
\author{Igor Prokhorenkov}
\author{Ken Richardson}
\address{Department of Mathematics\\
Texas Christian University\\
Box 298900 \\
Fort Worth, Texas 76129}
\email{i.prokhorenkov@tcu.edu\\
k.richardson@tcu.edu}
\date{May, 2008}
\maketitle

\begin{abstract}
We introduce a new class of natural, explicitly defined, transversally
elliptic differential operators over manifolds with compact group actions.
Under certain assumptions, the symbols of these operators generate all the
possible values of the equivariant index. We also show that the components of the
representation-valued equivariant index coincide with those of an elliptic
operator constructed from the original data.
\end{abstract}

\tableofcontents

\section{Introduction}

\vspace{0in}The representation-valued equivariant index of a transversally
elliptic operator is an important invariant in $K$-theory (see \cite{A}).
There are few known nontrivial examples in the literature where this
invariant is explicitly computed. Part of the motivation of this paper is to
provide an interesting and sufficiently general class of examples of
transversally elliptic differential operators for which such computations
are possible.

It is well-known that each compactly supported $K$-theory class of the
cotangent bundle over an even-dimensional spin$^{c}$ manifold is represented
by the symbol of a Dirac-type operator. This implies that Dirac-type symbols
map onto the image of the $K$-theory index homomorphism. In this paper, we
will generalize the second fact to the case of transversally elliptic
operators over a compact manifold endowed with a compact Lie group action.
The role of the Dirac-type operators will be played by a new class of
transversally elliptic differential operators introduced in this paper. To
construct these operators, we lift the group action to a principal bundle so
that all orbits in the principal bundle have the same dimension. There is a
natural transversal Dirac operator associated to this action. This operator
induces a transversally elliptic differential operator on the base manifold
with the desired properties. In the case when all orbits have the same
dimension, the orbits on the base manifold form a Riemannian foliation, and
our construction produces a transversal Dirac operator as studied by \cite%
{GlK}, \cite{Hab}, \cite{BrKRi}, \cite{DouglasGlK}, \cite{Jung1}, \cite%
{Jung2}, and others. This new operator will generate all possible values of
the representation-valued equivariant index. Further, we show that the
decomposition of this equivariant index representation into irreducible
components may be computed by means of equivariant indices of elliptic
operators. Thus, the techniques of Atiyah and Segal \cite{ASe} for elliptic
operators become applicable to transversally elliptic operators as well.

Now we describe the content of the paper. Let $M$ be a closed Riemannian
manifold. Let $Q\subset TM$ be a smooth distribution over $M$; we do not
assume that $Q$ or its normal bundle are involutive. Section \ref%
{RestrCliffBundlesSection} contains preliminary results about connections
associated to restrictions of Clifford structures. In Section \ref%
{transverseDiracDistributions}, we assume only that $E\rightarrow M$ is a $%
\mathbb{C}\mathrm{l}\left( Q\right) $-bundle with corresponding compatible
Clifford connection. Such a connection always exists if $E$ is in addition a 
$\mathbb{C}\mathrm{l}\left( TM\right) $-bundle; see Section \ref%
{RestrCliffBundlesSection}. Using this $\mathbb{C}\mathrm{l}\left( Q\right) $
connection, we construct an operator $D_{Q}$ whose principal symbol $\sigma
\left( D_{Q}\right) \left( \xi _{x}\right) $ is invertible for all $\xi
_{x}\in Q_{x}\setminus \left\{ 0\right\} $ and prove that it is essentially
self-adjoint. In the case of a Riemannian foliation with normal bundle $Q$,
this construction produces the well-known self-adjoint version of the
transversal Dirac operator (see \cite{GlK}, \cite{Jung1}, etc.).

In Section \ref{EquivariantOpsFrameBndleSection}, we assume that there is an
isometric action of a compact Lie group $G$ on $M$. The action of $G$ lifts
to the orthonormal frame bundle $F_{O}$ of $M$. Given an equivariant
transversal Dirac operator on $F_{O}$ and irreducible representation of the
orthogonal group, we show how to construct a transversally elliptic operator
on $M$. Similarly, using the transversal Dirac operator and an irreducible
representation of $G$, we construct an elliptic operator on $F_{O}\diagup G$%
. The precise relationship between the eigenspaces of these two operators is
stated in Proposition \ref{SpectraSame}.

In Section \ref{topPropertiesSection}, we study the equivariant index of the 
$K$-theory class of operators constructed in Section \ref%
{EquivariantOpsFrameBndleSection}. To do this, we first derive a
multiplicative property of the equivariant index on associated fiber bundles
with compact fibers. This property, stated in Theorem \ref%
{EquivariantMultiplicativeTheorem}, is a generalization of the
multiplicative property of the index for sphere bundles shown by Atiyah and
Singer in \cite{ASi1}. The main result of Section \ref{topPropertiesSection}
is Theorem \ref{indexClassGivenByTransvDiracThm}, in which we show that the
symbols of the lifted transversal Dirac operators generate all the possible
equivariant indices, if the $F_{O}$ is $G$-transversally spin$^{c}$.

In Section \ref{ExampleSection}, we demonstrate our constructions of the
lifted transversal Dirac operator and verify our results by explicit
calculations on the two-sphere.

The reader may consult \cite{A} for the basic properties of transversally
elliptic equivariant operators and their equivariant indices. Interesting
and relevant results also appear in \cite{ASe}, \cite{B-G-V}, \cite{Be-V1}, 
\cite{Be-V2}, \cite{Bra3}, \cite{BrKRi}, and \cite{Paradan}.

\section{Restrictions of Clifford structures\label{RestrCliffBundlesSection}}

Let $M$ be a closed Riemannian manifold with metric $\left\langle ~\cdot
~,~\cdot ~\right\rangle $, and let $E$ be a Clifford bundle over $M$. Recall
that a Clifford bundle $E$ is a complex Hermitian vector bundle endowed with
a Clifford action $c:TM\otimes \mathbb{C}\rightarrow \mathrm{End}\left(
E\right) $ and a connection $\nabla ^{E}$ compatible with this action and
the metric. Let $Q$ be a subbundle of $TM$, and let $L$ be the orthogonal
complement of $Q$ in $TM$. The Levi-Civita connection $\nabla ^{M}$ induces
a connection $\nabla ^{Q}$ on $Q$ by the following formula. Given any
section $Y\in \Gamma Q$ and any vector field $X\in \Gamma \left( TM\right) $%
, define%
\begin{equation}
\nabla _{X}^{Q}Y=\pi \nabla _{X}^{M}Y,  \label{Qconnection}
\end{equation}%
where $\pi :TM\rightarrow Q$ is the orthogonal bundle projection. It is
elementary to check that Formula (\ref{Qconnection}) yields a metric
connection on $Q$ (with the restricted metric).

\begin{remark}
Note that this connection is not generally $Q$-torsion-free, because the
torsion-free property is equivalent to the integrability of $Q$ (i.e. $\left[
\Gamma Q,\Gamma Q\right] \subset \Gamma Q$).
\end{remark}

We now modify the connection $\nabla ^{E}$ so that it has the desired
compatibility with $\nabla ^{Q}$. Every metric connection $\widetilde{\nabla
^{E}}$ on $E$ satisfies%
\begin{equation*}
\widetilde{\nabla _{X}^{E}}=\nabla _{X}^{E}+B_{X},
\end{equation*}%
where $B_{X}$ is a skew-Hermitian endomorphism of $E$ that is $C^{\infty
}\left( M\right) $-linear in $X$. In order that $\widetilde{\nabla ^{E}}$ is
a $\mathbb{C}\mathrm{l}\left( Q\right) $-connection compatible with $\nabla
^{Q}$, we must have that if $X\in \Gamma \left( TM\right) $, $Y\in \Gamma Q$%
, $s\in \Gamma E$,%
\begin{equation*}
\widetilde{\nabla _{X}^{E}}c\left( Y\right) s=c\left( \nabla
_{X}^{Q}Y\right) s+c\left( Y\right) \widetilde{\nabla _{X}^{E}}s\text{.}
\end{equation*}%
By (\ref{Qconnection}) we see that%
\begin{equation*}
B_{X}c\left( Y\right) s=-c\left( \left( 1-\pi \right) \nabla
_{X}^{M}Y\right) s+c\left( Y\right) B_{X}s,
\end{equation*}%
or%
\begin{equation}
c\left( \left( 1-\pi \right) \nabla _{X}^{M}Y\right) s=\left[ c\left(
Y\right) ,B_{X}\right] s.  \label{cliffordRequirement}
\end{equation}%
Computing with a local orthonormal frame $\left\{ e_{1},...,e_{p}\right\} $
for $L$, we have 
\begin{eqnarray*}
\left( 1-\pi \right) \nabla _{X}^{M}Y &=&\sum_{m=1}^{p}\left\langle \nabla
_{X}^{M}Y,e_{m}\right\rangle e_{m}=-\sum_{m=1}^{p}\left\langle Y,\nabla
_{X}^{M}e_{m}\right\rangle e_{m} \\
&=&-\sum_{m=1}^{p}\left\langle Y,\pi \nabla _{X}^{M}e_{m}\right\rangle e_{m}
\\
&=&\frac{1}{2}\sum_{m=1}^{p}\left( c\left( Y\right) c\left( \pi \nabla
_{X}^{M}e_{m}\right) +c\left( \pi \nabla _{X}^{M}e_{m}\right) c\left(
Y\right) \right) e_{m}~.
\end{eqnarray*}%
Then (\ref{cliffordRequirement}) implies that%
\begin{equation*}
\frac{1}{2}\sum_{m=1}^{p}\left( c\left( Y\right) c\left( \pi \nabla
_{X}^{M}e_{m}\right) c\left( e_{m}\right) -c\left( \pi \nabla
_{X}^{M}e_{m}\right) c\left( e_{m}\right) c\left( Y\right) \right) =\left(
c\left( Y\right) B_{X}-B_{X}c\left( Y\right) \right) ,
\end{equation*}%
or%
\begin{equation*}
\left( B_{X}-\frac{1}{2}\sum_{m=1}^{p}c\left( \pi \nabla
_{X}^{M}e_{m}\right) c\left( e_{m}\right) \right) c\left( Y\right) =c\left(
Y\right) \left( B_{X}-\frac{1}{2}\sum_{m=1}^{p}c\left( \pi \nabla
_{X}^{M}e_{m}\right) c\left( e_{m}\right) \right) .
\end{equation*}%
We conclude that the requirement that $\widetilde{\nabla ^{E}}$ is a $%
\mathbb{C}\mathrm{l}\left( Q\right) $-connection defines $B_{X}$ (and thus $%
\widetilde{\nabla _{X}^{E}}$ ) up to a skew-adjoint endomorphism of $E$ that
commutes with Clifford multiplication by vectors in $Q$. We may always take%
\begin{equation*}
B_{X}=\frac{1}{2}\sum_{m=1}^{p}c\left( \pi \nabla _{X}^{M}e_{m}\right)
c\left( e_{m}\right) .
\end{equation*}%
This choice of $B_{X}$ (and thus $\widetilde{\nabla _{X}^{E}}$ ) is
well-defined and canonical, since the formula is independent of the local
orthonormal frame for $L$.

We now show how to express $B_{X}$ in terms of any local orthonormal frame $%
f_{1},...,f_{q}$ for $Q$. 
\begin{eqnarray*}
B_{X} &=&\frac{1}{2}\sum_{m=1}^{p}c\left( \pi \nabla _{X}^{M}e_{m}\right)
c\left( e_{m}\right) =\frac{1}{2}\sum_{m=1}^{p}\sum_{j=1}^{q}\left\langle
\nabla _{X}^{M}e_{m},f_{j}\right\rangle c\left( f_{j}\right) c\left(
e_{m}\right) \\
&=&-\frac{1}{2}\sum_{m=1}^{p}\sum_{j=1}^{q}\left\langle e_{m},\nabla
_{X}^{M}f_{j}\right\rangle c\left( f_{j}\right) c\left( e_{m}\right) \\
&=&\frac{1}{4}\sum_{m=1}^{p}\sum_{j=1}^{q}\left\{ c\left( e_{m}\right)
c\left( \nabla _{X}^{M}f_{j}\right) +c\left( \nabla _{X}^{M}f_{j}\right)
c\left( e_{m}\right) \right\} c\left( f_{j}\right) c\left( e_{m}\right) \\
&=&\frac{1}{4}\sum_{j=1}^{q}\sum_{m=1}^{p}\left\{ -c\left( e_{m}\right)
c\left( \nabla _{X}^{M}f_{j}\right) c\left( e_{m}\right) c\left(
f_{j}\right) +c\left( \nabla _{X}^{M}f_{j}\right) c\left( f_{j}\right)
\right\} \\
&=&\frac{1}{4}\sum_{j=1}^{q}\left\{ -pc\left( \pi \nabla
_{X}^{M}f_{j}\right) c\left( f_{j}\right) -\left( p-2\right) c\left( \left(
1-\pi \right) \nabla _{X}^{M}f_{j}\right) c\left( f_{j}\right) +pc\left(
\nabla _{X}^{M}f_{j}\right) c\left( f_{j}\right) \right\} \\
&=&\frac{1}{2}\sum_{j=1}^{q}c\left( \left( 1-\pi \right) \nabla
_{X}^{M}f_{j}\right) c\left( f_{j}\right) .
\end{eqnarray*}%
Observe that this expression for $B_{X}$ is the same as the original
expression for $B_{X}$ with $Q$ replaced by $L$. We have shown the following.

\begin{proposition}
Let $M$ be a closed Riemannian manifold, and let $\left( E,\nabla
^{E},c\right) $ be a Hermitian Clifford bundle over $M$. Let $Q$ be a
subbundle of $TM$, and let $\nabla ^{Q}$ denote the metric connection on $Q$
defined by $\nabla _{X}^{Q}Y=\pi \nabla _{X}^{M}Y$, where $\pi
:TM\rightarrow Q$ is the orthogonal bundle projection. Then the connection $%
\widetilde{\nabla ^{E}}$ defined by%
\begin{equation*}
\widetilde{\nabla _{X}^{E}}=\nabla _{X}^{E}+\frac{1}{2}\sum_{m=1}^{p}c\left(
\pi \nabla _{X}^{M}e_{m}\right) c\left( e_{m}\right)
\end{equation*}%
for all $X\in \Gamma \left( TM\right) $ is a well-defined $\mathbb{C}\mathrm{%
l}\left( Q\right) $-connection and a metric connection on $E$ with respect
to the connection $\nabla ^{Q}$. Furthermore, $\widetilde{\nabla ^{E}}$ is a 
$\mathbb{C}\mathrm{l}\left( L\right) $-connection and a metric connection on 
$E$ with respect to the connection $\nabla ^{L}=\left( 1-\pi \right) \nabla
^{M}$.
\end{proposition}

\section{Transverse Dirac Operators for Distributions\label%
{transverseDiracDistributions}}

We showed in Section \ref{RestrCliffBundlesSection} that, for a given
distribution $Q\subset TM$, it is always possible to obtain a bundle of $%
\mathbb{C}\mathrm{l}\left( Q\right) $-modules with Clifford connection from
a bundle of $\mathbb{C}\mathrm{l}\left( TM\right) $-Clifford modules. In
this section, we will assume more generally that a $\mathbb{C}\mathrm{l}%
\left( Q\right) $-module structure on a complex Hermitian vector bundle $E$
is given and will define transverse Dirac operators on sections of $E$. As
in Section \ref{RestrCliffBundlesSection}, $M$ is a closed Riemannian
manifold with metric $\left\langle ~\cdot ~,~\cdot ~\right\rangle $, $%
c:Q\rightarrow \mathrm{End}\left( E\right) $ is the Clifford multiplication
on $E$, and $\nabla ^{E}$ is a $\mathbb{C}\mathrm{l}\left( Q\right) $
connection that is compatible with the metric on $M$; that is, Clifford
multiplication by each vector is skew-Hermitian, and we require%
\begin{equation*}
\nabla _{X}^{E}\left( c\left( V\right) s\right) =c\left( \nabla
_{X}^{Q}V\right) s+c\left( V\right) \nabla _{X}^{E}s
\end{equation*}%
for all $X\in \Gamma \left( TM\right) $, $V\in \Gamma Q$, and $s\in \Gamma E$%
. Note that the connection $\widetilde{\nabla ^{E}}$ from Section \ref%
{RestrCliffBundlesSection} is an example of such a connection, but not all
such $\mathbb{C}\mathrm{l}\left( Q\right) $ connections are of that type.
Let $L=Q^{\bot }$, let $\left( f_{1},...,f_{q}\right) $ be a local
orthonormal frame for $Q$, and let $\pi :TM\rightarrow Q$ be the orthogonal
projection. We define the Dirac operator $A_{Q}$ corresponding to the
distribution $Q$ as 
\begin{equation}
A_{Q}=\sum_{j=1}^{q}c\left( f_{j}\right) \nabla _{f_{j}}^{E}.  \label{AQdef}
\end{equation}%
This definition is independent of the choices made; in fact it is the
composition of the maps 
\begin{equation*}
\Gamma \left( E\right) \overset{\nabla ^{E}}{\rightarrow }\Gamma \left(
T^{\ast }M\otimes E\right) \overset{\cong }{\rightarrow }\Gamma \left(
TM\otimes E\right) \overset{\pi }{\rightarrow }\Gamma \left( Q\otimes
E\right) \overset{c}{\rightarrow }\Gamma \left( E\right) .
\end{equation*}%
We calculate the formal adjoint. Letting $\left( s_{1},s_{2}\right) $ denote
the pointwise inner product of sections of $E$, we have that 
\begin{eqnarray*}
\left( A_{Q}s_{1},s_{2}\right) -\left( s_{1},A_{Q}s_{2}\right)
&=&\sum_{j=1}^{q}\left( c\left( f_{j}\right) \nabla
_{f_{j}}^{E}s_{1},s_{2}\right) -\left( s_{1},c\left( f_{j}\right) \nabla
_{f_{j}}^{E}s_{2}\right) \\
&=&\sum_{j=1}^{q}\left( \nabla _{f_{j}}^{E}\left( c\left( f_{j}\right)
s_{1}\right) ,s_{2}\right) -\left( c\left( \pi \nabla
_{f_{j}}^{M}f_{j}\right) s_{1},s_{2}\right) +\left( c\left( f_{j}\right)
s_{1},\nabla _{f_{j}}^{E}s_{2}\right) \\
&=&\sum_{j=1}^{q}\nabla _{f_{j}}^{M}\left( c\left( f_{j}\right)
s_{1},s_{2}\right) -\left( c\left( \sum_{j=1}^{q}\pi \nabla
_{f_{j}}^{M}f_{j}\right) s_{1},s_{2}\right) \\
&=&-\sum_{j=1}^{q}\nabla _{f_{j}}^{M}i_{f_{j}}\omega +\omega \left(
\sum_{j=1}^{q}\pi \nabla _{f_{j}}^{M}f_{j}\right) ,
\end{eqnarray*}%
where $\omega $ is the one-form defined by $\omega \left( X\right) =-\left(
c\left( X\right) s_{1},s_{2}\right) $ for $X\in \Gamma Q$ and is zero for $%
X\in \Gamma L$. Continuing, 
\begin{eqnarray*}
\left( A_{Q}s_{1},s_{2}\right) -\left( s_{1},A_{Q}s_{2}\right)
&=&-\sum_{j=1}^{q}\left( i_{f_{j}}\nabla _{f_{j}}^{M}+i_{\nabla
_{f_{j}}^{M}f_{j}}\right) \omega +\omega \left( \sum_{j=1}^{q}\pi \nabla
_{f_{j}}^{M}f_{j}\right) \\
&=&-\sum_{j=1}^{q}\left( i_{f_{j}}\nabla _{f_{j}}^{M}+i_{\left( \pi \nabla
_{f_{j}}^{M}f_{j}\right) }\right) \omega +\omega \left( \sum_{j=1}^{q}\pi
\nabla _{f_{j}}^{M}f_{j}\right) \\
&=&-\sum_{j=1}^{q}i_{f_{j}}\nabla _{f_{j}}^{M}\omega
=-\sum_{j=1}^{q}i_{f_{j}}\pi \nabla _{f_{j}}^{M}\omega ,
\end{eqnarray*}%
where the orthogonal projection $T^{\ast }M\rightarrow Q^{\ast }$ is denoted
by $\pi $ as well. In what follows, let $\left( e_{1},...,e_{p}\right) $ be
an orthonormal frame of $L$, and let $\nabla ^{M}=$ $\nabla ^{Q}+\nabla
^{L}=\pi \nabla ^{M}+\left( 1-\pi \right) \nabla ^{M}$ on forms. The
divergence of a general one-form $\beta $ that is zero on $L$ is 
\begin{eqnarray*}
\delta \beta &=&-\sum_{j=1}^{q}i_{f_{j}}\nabla _{f_{j}}^{M}\beta
-\sum_{m=1}^{p}i_{e_{m}}\nabla _{e_{m}}^{M}\beta \\
&=&-\sum_{j=1}^{q}i_{f_{j}}\pi \nabla _{f_{j}}^{M}\beta
-\sum_{m=1}^{p}i_{e_{m}}\nabla _{e_{m}}^{L}\beta ,
\end{eqnarray*}%
Letting $\beta =\sum_{k=1}^{q}\beta _{k}f_{k}^{\ast },$ then 
\begin{eqnarray*}
\delta \beta +\sum_{j=1}^{q}i_{f_{j}}\pi \nabla _{f_{j}}^{M}\beta
&=&-\sum_{m=1}^{p}i_{e_{m}}\nabla _{e_{m}}^{L}\left( \sum_{k=1}^{q}\beta
_{k}f_{k}^{\ast }\right) \\
&=&-\sum_{k=1}^{q}\sum_{m=1}^{p}\beta _{k}i_{e_{m}}\nabla _{e_{m}}^{L}\left(
f_{k}^{\ast }\right) \\
&=&-\sum_{k=1}^{q}\sum_{m=1}^{p}\beta _{k}i_{e_{m}}\left(
\sum_{j=1}^{p}\left( \nabla _{e_{m}}^{M}\left( f_{k}^{\ast }\right)
,e_{j}^{\ast }\right) e_{j}^{\ast }\right) \\
&=&\sum_{k=1}^{q}\sum_{m=1}^{p}\beta _{k}i_{e_{m}}\left(
\sum_{j=1}^{p}\left( \nabla _{e_{m}}^{M}\left( e_{j}^{\ast }\right)
,f_{k}^{\ast }\right) e_{j}^{\ast }\right) \\
&=&\sum_{k=1}^{q}\beta _{k}\left( \sum_{m=1}^{p}\left( \nabla
_{e_{m}}^{M}\left( e_{m}^{\ast }\right) ,f_{k}^{\ast }\right) \right) \\
&=&i_{H^{L}}\beta ,
\end{eqnarray*}%
where $H^{L}$ is the mean curvature vector field of $L$. Thus, for every
one-form $\beta $ that is zero on $L$ , 
\begin{equation*}
-\sum_{j=1}^{q}i_{f_{j}}\pi \nabla _{f_{j}}^{M}\beta =\delta \beta
-i_{H^{L}}\beta .
\end{equation*}%
Applying this result to the form $\omega $ defined above, we have 
\begin{eqnarray*}
\left( A_{Q}s_{1},s_{2}\right) -\left( s_{1},A_{Q}s_{2}\right) &=&\delta
\omega -i_{H^{L}}\omega \\
&=&\delta \omega -\left( s_{1},c\left( H^{L}\right) s_{2}\right) .
\end{eqnarray*}%
Thus, the formal adjoint $A_{Q}^{\ast }$ of $A_{Q}$ is%
\begin{equation*}
A_{Q}^{\ast }=A_{Q}-c\left( H^{L}\right) ,
\end{equation*}%
and the operator%
\begin{equation}
D_{Q}=A_{Q}-\frac{1}{2}c\left( H^{L}\right)  \label{DQdef}
\end{equation}%
is formally self-adjoint.

A quick look at \cite{C} yields the following.

\begin{theorem}
For each distribution $Q\subset TM$ and every bundle $E$ of $\mathbb{C}%
\mathrm{l}\left( Q\right) $-modules, the transversally elliptic operator $%
D_{Q}$ defined by (\ref{AQdef}) and (\ref{DQdef}) is essentially
self-adjoint.
\end{theorem}

It is not necessarily the case that the spectrum of $D_{Q}$ is discrete, as
the following example shows.

\begin{example}
We consider the torus $M=\left( \mathbb{R}\diagup 2\pi \mathbb{Z}\right)
^{2} $ with the metric $e^{2g\left( y\right) }dx^{2}+dy^{2}$ for some $2\pi $%
-periodic smooth function $g$. Consider the orthogonal distributions $L=%
\mathrm{span}\left\{ \partial _{y}\right\} $ and $Q=\mathrm{span}\left\{
\partial _{x}\right\} $. Let $E$ be the trivial complex line bundle over $M$%
, and let $\mathbb{C}\mathrm{l}\left( Q\right) $ and $\mathbb{C}\mathrm{l}%
\left( L\right) $ both act on $E$ via $c\left( \partial _{y}\right)
=i=c\left( e^{-g\left( y\right) }\partial _{x}\right) $. The connections $%
\nabla ^{L}$ and $\nabla ^{Q}$ satisfy 
\begin{eqnarray*}
\nabla _{\partial _{y}}^{L}\partial _{y} &=&\nabla _{\partial
_{x}}^{L}\partial _{y}=0; \\
\nabla _{\partial _{x}}^{Q}\partial _{x} &=&0;~\nabla _{\partial
_{y}}^{Q}\partial _{x}=g^{\prime }\left( y\right) \partial _{x}.
\end{eqnarray*}

The trivial connection $\nabla ^{E}$ is a $\mathbb{C}\mathrm{l}\left(
L\right) $ connection with respect to $\nabla ^{L}$ and is also a $\mathbb{C}%
\mathrm{l}\left( Q\right) $ connection with respect to $\nabla ^{Q}$.
Observe that the mean curvatures of these distributions are%
\begin{eqnarray*}
H^{Q} &=&\left( 1-\pi \right) \nabla _{e^{-g\left( y\right) }\partial
_{x}}^{M}e^{-g\left( y\right) }\partial _{x}=-g^{\prime }\left( y\right)
\partial _{y}\text{ and} \\
H^{L} &=&\pi \nabla _{\partial _{y}}^{M}\partial _{y}=\nabla _{\partial
_{y}}^{M}\partial _{y}=0
\end{eqnarray*}%
From formulas (\ref{AQdef}) and (\ref{DQdef}),%
\begin{eqnarray*}
A_{L} &=&i\partial _{y},\text{ and} \\
D_{L} &=&i\left( \partial _{y}+\frac{1}{2}g^{\prime }\left( y\right) \right)
.
\end{eqnarray*}%
The spectrum $\sigma \left( D_{L}\right) =\mathbb{Z}$ is a set consisting of
eigenvalues of infinite multiplicity, and thus $\sigma \left( D_{L}\right) $
consists entirely of pure point spectrum. The eigenspace $E_{n}$
corresponding to the eigenvalue $n$ is 
\begin{equation*}
E_{n}=\left\{ e^{-iny-\frac{g\left( y\right) }{2}}f\left( x\right) :f\in
L^{2}\left( S^{1}\right) \right\} ,
\end{equation*}%
and $\bigcup\limits_{n\in \mathbb{Z}}E_{n}$ is dense in $L^{2}\left(
M\right) $.

On the other hand, the operator%
\begin{equation*}
D_{Q}=A_{Q}-\frac{1}{2}c\left( H^{L}\right) =A_{Q}=ie^{-g\left( y\right)
}\partial _{x}
\end{equation*}%
has only one eigenvalue, $0$, corresponding to the eigenspace $\left\{
h\left( y\right) :h\in L^{2}\left( S^{1}\right) \right\} $. Next, note that $%
F_{n}=\left\{ e^{-inx}\psi \left( y\right) :\psi \in L^{2}\left(
S^{1}\right) \right\} $ is an invariant subspace for $D_{Q}$ . The spectrum
of the restriction of $D_{Q}$ to $F_{n}$ is $n\left[ a,b\right] $, where $%
\left[ a,b\right] \subset \left( 0,\infty \right) $ is the range of $%
e^{-g\left( y\right) }$. Thus, the spectrum $\sigma \left( D_{Q}\right) $ is%
\begin{equation*}
\sigma \left( D_{Q}\right) =\bigcup\limits_{n\in \mathbb{Z}}n\left[ a,b%
\right] ,
\end{equation*}%
and the pure point spectrum of $D_{Q}$ is $\left\{ 0\right\} $.
\end{example}

\begin{example}
\label{RiemFoliationTrDiracExample}Suppose a closed manifold $M$ is endowed
with a Riemannian foliation $\mathcal{F}$ such that the metric is
bundlelike, meaning that the leaves are locally equidistant. If the orbits
of a $G$-manifold have the same dimension, then they form a Riemannian
foliation. In such foliations, there is a natural construction of
transversal Dirac operators (see \cite{BrKRi} , \cite{GlK} , \cite{La}),
which is a special case of the construction in this section. Choose a local
adapted frame field $\left\{ e_{1},...,e_{n}\right\} $ for the tangent
bundle of $M$ , such that $\left\{ e_{1},...,e_{q}\right\} $ is a local
basis of the normal bundle $N\mathcal{F}$ for the foliation and such that
each $e_{j}$ is a basic vector field for $1\leq j\leq q$. The word \emph{%
basic} means that the flows of those vector fields map leaves to leaves, and
such a basis can be chosen near every point if and only if the foliation is
Riemannian. Next, assume that we have a complex Hermitian vector bundle $%
E\rightarrow M$ that is a bundle of $\mathbb{C}\mathrm{l}\left( N\mathcal{F}%
\right) $ modules that is equivariant with respect to the $G$ action, and
let $\nabla $ be the corresponding equivariant, metric, Clifford connection.
We define the transversal Dirac operator by 
\begin{equation*}
A_{N\mathcal{F}}=\sum_{j=1}^{q}c\left( e_{j}\right) \nabla _{e_{j}},
\end{equation*}%
as in the notation of this section. As before, the operator 
\begin{equation*}
D_{N\mathcal{F}}=A_{N\mathcal{F}}-\frac{1}{2}c\left( H\right)
\end{equation*}%
is an essentially self-adjoint operator, where $H$ is the mean curvature
vector field of the orbits.
\end{example}

\section{Equivariant operators on the frame bundle\label%
{EquivariantOpsFrameBndleSection}}

\subsection{Equivariant structure of the orthonormal frame bundle\label%
{equivariantStructureSection}}

\vspace{0in}Given a complete, connected $G$-manifold, the action of $g\in G$
on $M$ induces an action of $dg$ on $TM$, which in turn induces an action of 
$G$ on the principal $O\left( n\right) $-bundle $F_{O}\overset{p}{%
\rightarrow }M$ of orthonormal frames over $M$.

\begin{lemma}
\vspace{0in}The action of $G$ on $F_{O}$ is regular, i.e. the isotropy
subgroups corresponding to any two points of $M$ are conjugate.
\end{lemma}

\begin{proof}
Let $H$ be the isotropy subgroup of a frame $f\in F_{O}$. Then $H$ also
fixes $p\left( f\right) \in M$, and since $H$ fixes the frame, its
differentials fix the entire tangent space at $p\left( f\right) $. Since it
fixes the tangent space, every element of $H$ also fixes every frame in $%
p^{-1}\left( p\left( f\right) \right) $; thus every frame in a given fiber
must have the same isotropy subgroup. Since the elements of $H$ map
geodesics to geodesics and preserve distance, a neighborhood of $p\left(
f\right) $ is fixed by $H$. Thus, $H$ is a subgroup of the isotropy subgroup
at each point of that neighborhood. Conversely, if an element of $G$ fixes a
neighborhood of a point $x$ in $M$, then it fixes all frames in $%
p^{-1}\left( x\right) $, and thus all frames in the fibers above that
neighborhood. Since $M$ is connected, we may conclude that every point of $%
F_{O}$ has the same isotropy subgroup $H$, and $H$ is the subgroup of $G$
that fixes every point of $M$.
\end{proof}

\begin{remark}
Since this subgroup $H$ is normal, we often reduce the group $G$ to the
group $G/H$ so that our action is effective, in which case the isotropy
subgroups on $F_{O}$ are all trivial.
\end{remark}

In any case, the $G$ orbits on $F_{O}$ are diffeomorphic and form a
Riemannian fiber bundle, in the natural metric on $F_{O}$ defined as
follows. The Levi-Civita connection on $M$ determines the horizontal
subbundle $\mathcal{H}$ of $TF_{O}$. We construct the local product metric
on $F_{O}$ using a biinvariant fiber metric and the pullback of the metric
on $M$ to $\mathcal{H}$; with this metric, $F_{O}$ is a compact Riemannian $%
G\times O\left( n\right) $-manifold. The lifted $G$-action commutes with the 
$O\left( n\right) $-action. Let $\mathcal{F}$ denote the foliation of $G$%
-orbits on $F_{O}$, and observe that $F_{O}\overset{\pi }{\rightarrow }%
F_{O}\diagup G=F_{O}\diagup \mathcal{F}$ is a Riemannian submersion of
compact $O\left( n\right) $-manifolds.

Let $E\rightarrow F_{O}$ be a Hermitian vector bundle that is equivariant
with respect to the $G\times O\left( n\right) $ action. Let $\rho
:G\rightarrow U\left( V_{\rho }\right) $ and $\sigma :O\left( n\right)
\rightarrow U\left( W_{\sigma }\right) $ be irreducible unitary
representations. We define the bundle $\mathcal{E}^{\sigma }\rightarrow M$
by 
\begin{equation*}
\mathcal{E}_{x}^{\sigma }=\Gamma \left( p^{-1}\left( x\right) ,E\right)
^{\sigma },
\end{equation*}%
where the superscript $\sigma $ is defined for a $O\left( n\right) $-module $%
Z$ by%
\begin{equation*}
Z^{\sigma }=\mathrm{eval}\left( \mathrm{Hom}_{O\left( n\right) }\left(
W_{\sigma },Z\right) \otimes W_{\sigma }\right) ,
\end{equation*}%
where $\mathbb{\mathrm{eval}}:\mathrm{Hom}_{O\left( n\right) }\left(
W_{\sigma },Z\right) \otimes W_{\sigma }\rightarrow Z$ is the evaluation map 
$\phi \otimes w\mapsto \phi \left( w\right) $. The space $Z^{\sigma }$ is
the vector subspace of $Z$ on which $O\left( n\right) $ acts as a direct sum
of representations of type $\sigma $. The bundle $\mathcal{E}^{\sigma }$ is
a Hermitian $G$-vector bundle of finite rank over $M$. The metric on $%
\mathcal{E}^{\sigma }$ is chosen as follows. For any $v_{x}$,$w_{x}\in 
\mathcal{E}_{x}^{\sigma }$, we define%
\begin{equation*}
\,\left\langle v_{x},w_{x}\right\rangle :=\int_{p^{-1}\left( x\right)
}\left\langle v_{x}\left( y\right) ,w_{x}\left( y\right) \right\rangle
_{y,E}~d\mu _{x}\left( y\right) ,
\end{equation*}%
where $d\mu _{x}$ is the measure on $p^{-1}\left( x\right) $ induced from
the metric on $F_{O}$. See \cite{BrKRi} for a similar construction.

Similarly, we define the bundle $\mathcal{T}^{\rho }\rightarrow F_{O}\diagup
G$ by%
\begin{equation*}
\mathcal{T}_{y}^{\rho }=\Gamma \left( \pi ^{-1}\left( y\right) ,E\right)
^{\rho },
\end{equation*}%
and $\mathcal{T}^{\rho }\rightarrow F_{O}\diagup G$ is a Hermitian $O\left(
n\right) $-equivariant bundle of finite rank. The metric on $\mathcal{T}%
^{\rho }$ is%
\begin{equation*}
\left\langle v_{z},w_{z}\right\rangle :=\int_{\pi ^{-1}\left( y\right)
}\left\langle v_{z}\left( y\right) ,w_{z}\left( y\right) \right\rangle
_{z,E}~dm_{z}\left( y\right) ,
\end{equation*}%
where $dm_{z}$ is the measure on $\pi ^{-1}\left( z\right) $ induced from
the metric on $F_{O}$.

The vector spaces of sections $\Gamma \left( M,\mathcal{E}^{\sigma }\right) $
and $\Gamma \left( F_{O},E\right) ^{\sigma }$ can be identified via the
isomorphism%
\begin{equation*}
i_{\sigma }:\Gamma \left( M,\mathcal{E}^{\sigma }\right) \rightarrow \Gamma
\left( F_{O},E\right) ^{\sigma },
\end{equation*}%
where for any section $s\in \Gamma \left( M,\mathcal{E}^{\sigma }\right) $, $%
s\left( x\right) \in \Gamma \left( p^{-1}\left( x\right) ,E\right) ^{\sigma
} $ for each $x\in M$, and we let%
\begin{equation*}
i_{\sigma }\left( s\right) \left( f_{x}\right) :=\left. s\left( x\right)
\right\vert _{f_{x}}
\end{equation*}%
for every $f_{x}\in p^{-1}\left( x\right) \subset F_{O}$. Then $i_{\sigma
}^{-1}:\Gamma \left( F_{O},E\right) ^{\sigma }\rightarrow \Gamma \left( M,%
\mathcal{E}^{\sigma }\right) $ is given by%
\begin{equation*}
i_{\sigma }^{-1}\left( u\right) \left( x\right) =\left. u\right\vert
_{p^{-1}\left( x\right) }.
\end{equation*}%
Observe that $i_{\sigma }:\Gamma \left( M,\mathcal{E}^{\sigma }\right)
\rightarrow \Gamma \left( F_{O},E\right) ^{\sigma }$ extends to an $L^{2}$
isometry. Given $u,v\in $ $\Gamma \left( M,\mathcal{E}^{\sigma }\right) $,%
\begin{eqnarray*}
\left\langle u,v\right\rangle _{M} &=&\int_{M}\left\langle
u_{x},v_{x}\right\rangle ~dx=\int_{M}\int_{p^{-1}\left( x\right)
}\left\langle u_{x}\left( y\right) ,v_{x}\left( y\right) \right\rangle
_{y,E}~d\mu _{x}\left( y\right) ~dx \\
&=&\int_{M}\left( \int_{p^{-1}\left( x\right) }\left\langle i_{\sigma
}\left( u\right) ,i_{\sigma }\left( v\right) \right\rangle _{E}~d\mu
_{x}\left( y\right) \right) ~dx \\
&=&\int_{F_{O}}\left\langle i_{\sigma }\left( u\right) ,i_{\sigma }\left(
v\right) \right\rangle _{E}~=\left\langle i_{\sigma }\left( u\right)
,i_{\sigma }\left( v\right) \right\rangle _{F_{O}},
\end{eqnarray*}%
where $dx$ is the Riemannian measure on $M$; we have used the fact that $p$
is a Riemannian submersion. Similarly, we let%
\begin{equation*}
j_{\rho }:\Gamma \left( F_{O}\diagup G,\mathcal{T}^{\rho }\right)
\rightarrow \Gamma \left( F_{O},E\right) ^{\rho }
\end{equation*}%
be the natural identification, which extends to an $L^{2}$ isometry.

Let%
\begin{equation*}
\Gamma \left( M,\mathcal{E}^{\sigma }\right) ^{\alpha }=\mathrm{eval}\left( 
\mathrm{Hom}_{G}\left( V_{\alpha },\Gamma \left( M,\mathcal{E}^{\sigma
}\right) \right) \otimes V_{\alpha }\right) .
\end{equation*}%
Similarly, let%
\begin{equation*}
\Gamma \left( F_{O}\diagup G,\mathcal{T}^{\rho }\right) ^{\beta }=\mathrm{%
eval}\left( \mathrm{Hom}_{G}\left( W_{\beta },\Gamma \left( F_{O}\diagup G,%
\mathcal{T}^{\rho }\right) \right) \otimes W_{\beta }\right) .
\end{equation*}

\begin{theorem}
\label{IsomorphismsOfSectionsTheorem}For any irreducible representations $%
\rho :G\rightarrow U\left( V_{\rho }\right) $ and $\sigma :O\left( n\right)
\rightarrow U\left( W_{\sigma }\right) $, the map $j_{\rho }^{-1}\circ
i_{\sigma }:\Gamma \left( M,\mathcal{E}^{\sigma }\right) ^{\rho }\rightarrow
\Gamma \left( F_{O}\diagup G,\mathcal{T}^{\rho }\right) ^{\sigma }$ is an
isomorphism (with inverse $i_{\sigma }^{-1}\circ j_{\rho }$) that extends to
an $L^{2}$-isometry.
\end{theorem}

\begin{proof}
Observe that $i_{\sigma }$ implements the isomorphism 
\begin{eqnarray*}
\Gamma \left( M,\mathcal{E}^{\sigma }\right) &=&\Gamma \left( M,\Gamma
\left( p^{-1}\left( \cdot \right) ,\left. E\right\vert _{p^{-1}\left( \cdot
\right) }\right) ^{\sigma }\right) \\
&\cong &\Gamma \left( M,\Gamma \left( p^{-1}\left( \cdot \right) ,\left.
E\right\vert _{p^{-1}\left( \cdot \right) }\right) \right) ^{\sigma }=\Gamma
\left( F_{O},E\right) ^{\sigma }
\end{eqnarray*}%
to the space of sections of $E$ of $O\left( n\right) $ representation type $%
\sigma $. Its restriction to $\Gamma \left( M,\mathcal{E}^{\sigma }\right)
^{\rho }$ is%
\begin{eqnarray*}
\Gamma \left( M,\mathcal{E}^{\sigma }\right) ^{\rho } &=&\Gamma \left(
M,\Gamma \left( p^{-1}\left( \cdot \right) ,\left. E\right\vert
_{p^{-1}\left( \cdot \right) }\right) ^{\sigma }\right) ^{\rho } \\
&\cong &\left( \Gamma \left( M,\Gamma \left( p^{-1}\left( \cdot \right)
,\left. E\right\vert _{p^{-1}\left( \cdot \right) }\right) \right) ^{\sigma
}\right) ^{\rho } \\
&=&\left( \Gamma \left( F_{O},E\right) ^{\sigma }\right) ^{\rho }=\Gamma
\left( F_{O},E\right) ^{\sigma ,\rho },
\end{eqnarray*}%
where the superscript $\sigma ,\rho $ denotes restriction first to sections
of $O\left( n\right) $-representation type $\left[ \sigma \right] $ and then
to the subspace of sections of $G$-representation type $\left[ \rho \right] $%
. Since the $O\left( n\right) $ and $G$ actions commute, we may do this in
the other order, so that%
\begin{eqnarray*}
\Gamma \left( F_{O},E\right) ^{\sigma ,\rho } &=&\left( \Gamma \left(
F_{O},E\right) ^{\rho }\right) ^{\sigma } \\
&\cong &\Gamma \left( F_{O}\diagup G,\Gamma \left( \pi ^{-1}\left( y\right)
,\left. E\right\vert _{\pi ^{-1}\left( \cdot \right) }\right) ^{\rho
}\right) ^{\sigma } \\
&=&\Gamma \left( F_{O}\diagup G,\mathcal{T}^{\rho }\right) ^{\sigma },
\end{eqnarray*}%
where the isomorphism is the inverse of the restriction of $j_{\rho }$ to $%
\Gamma \left( F_{O}\diagup G,\mathcal{T}^{\rho }\right) ^{\sigma }$. Since $%
i_{\sigma }$ and $j_{\rho }$ are $L^{2}$ isometries, the result follows.
\end{proof}

\subsection{Dirac-type operators on the frame bundle\label{DiracFrameBundle}}

Let $E\rightarrow F_{O}$ be a Hermitian vector bundle of $\mathbb{C}\mathrm{l%
}\left( N\mathcal{F}\right) $ modules that is equivariant with respect to
the $G\times O\left( n\right) $ action. With notation as in Example \ref%
{RiemFoliationTrDiracExample}, we have the transversal Dirac operator $A_{N%
\mathcal{F}}$ defined by the composition 
\begin{equation*}
\Gamma \left( F_{O},E\right) \overset{\nabla }{\rightarrow }\Gamma \left(
F_{O},T^{\ast }F_{O}\otimes E\right) \overset{\mathrm{proj}}{\rightarrow }%
\Gamma \left( F_{O},N^{\ast }\mathcal{F}\otimes E\right) \overset{c}{%
\rightarrow }\Gamma \left( F_{O},E\right) .
\end{equation*}%
As explained previously, the operator 
\begin{equation*}
D_{N\mathcal{F}}=A_{N\mathcal{F}}-\frac{1}{2}c\left( H\right)
\end{equation*}%
is a essentially self-adjoint $G\times O\left( n\right) $-equivariant
operator, where $H$ is the mean curvature vector field of the $G$-orbits in $%
F_{O}$.

From $D_{N\mathcal{F}}$ we now construct equivariant differential operators
on $M$ and $F_{O}\diagup G$, as follows. We define the operators%
\begin{equation*}
D_{M}^{\sigma }:=i_{\sigma }^{-1}\circ D_{N\mathcal{F}}\circ i_{\sigma
}:\Gamma \left( M,\mathcal{E}^{\sigma }\right) \rightarrow \Gamma \left( M,%
\mathcal{E}^{\sigma }\right) ,
\end{equation*}%
and%
\begin{equation*}
D_{F_{O}\diagup G}^{\rho }:=j_{\rho }^{-1}\circ D_{N\mathcal{F}}\circ
j_{\rho }:\Gamma \left( F_{O}\diagup G,\mathcal{T}^{\rho }\right)
\rightarrow \Gamma \left( F_{O}\diagup G,\mathcal{T}^{\rho }\right) .
\end{equation*}%
For an irreducible representation $\alpha :G\rightarrow U\left( V_{\alpha
}\right) $, let 
\begin{equation*}
\left( D_{M}^{\sigma }\right) ^{\alpha }:\Gamma \left( M,\mathcal{E}^{\sigma
}\right) ^{\alpha }\rightarrow \Gamma \left( M,\mathcal{E}^{\sigma }\right)
^{\alpha }
\end{equation*}%
be the restriction of $D_{M}^{\sigma }$ to sections of $G$-representation
type $\left[ \alpha \right] $. Similarly, for an irreducible representation $%
\beta :G\rightarrow U\left( W_{\beta }\right) $, let 
\begin{equation*}
\left( D_{F_{O}\diagup G}^{\rho }\right) ^{\beta }:\Gamma \left(
F_{O}\diagup G,\mathcal{T}^{\rho }\right) ^{\beta }\rightarrow \Gamma \left(
F_{O}\diagup G,\mathcal{T}^{\rho }\right) ^{\beta }
\end{equation*}%
be the restriction of $D_{F_{O}\diagup G}^{\rho }$ to sections of $O\left(
n\right) $-representation type $\left[ \beta \right] $. The proposition
below follows from Theorem \ref{IsomorphismsOfSectionsTheorem}.

\begin{proposition}
\label{SpectraSame}The operator $D_{M}^{\sigma }$ is transversally elliptic
and $G$-equivariant, and $D_{F_{O}\diagup G}^{\rho }$ is elliptic and $%
O\left( n\right) $-equivariant, and the closures of these operators are
self-adjoint. The operators $\left( D_{M}^{\sigma }\right) ^{\rho }$ and $%
\left( D_{F_{O}\diagup G}^{\rho }\right) ^{\sigma }$ have identical discrete
spectrum, and the corresponding eigenspaces are conjugate via Hilbert space
isomorphisms.
\end{proposition}

Thus, questions about the transversally elliptic operator $D_{M}^{\sigma }$
can be reduced to questions about the elliptic operators $D_{F_{O}\diagup
G}^{\rho }$ for each irreducible $\rho :G\rightarrow U\left( V_{\rho
}\right) $.

\section{Topological properties of the lifted Dirac operators\label%
{topPropertiesSection}}

In this section, we will prove that if $F_{O}$ is $G$-transversally spin$^{c}
$, then the symbols of the lifted transversal Dirac operators generate all
the possible equivariant indices. To show this, we generalize the standard
multiplicative property of $K$-theory to the equivariant setting of our
paper.

\subsection{Equivariant Multiplicative Properties of $K$-theory\label%
{Multiplicative PropertySection}}

Let $H$ be a compact Lie group. Suppose that $P$ is a principal $H$-bundle
over a compact manifold $M$. Suppose that the compact Lie group $G$ acts on $%
M$ and lifts to $P$, such that the $G$-action on $P$ commutes with the $H$%
-action. Let $Z\overset{\pi }{\rightarrow }M$ be a fiber bundle associated
to $P$ with $H$-fiber $Y$; that is,%
\begin{equation*}
Z=P\times _{H}Y=P\times Y\diagup \left( p,y\right) \sim \left(
ph,h^{-1}y\right) .
\end{equation*}%
Then $G$ acts on $Z$ via $g\left[ \left( p,y\right) \right] =\left[ \left(
gp,y\right) \right] $.

For any $v\ $in the Lie algebra $\mathfrak{g}$ of $G$, let $\overline{v}$
denote the fundamental vector field on $M$ associated to $v$. As in \cite{A}%
, let 
\begin{equation*}
T_{G}^{\ast }M=\left\{ \xi \in T^{\ast }M:\xi \left( \overline{v}\right) =0%
\text{ for all }v\in \mathfrak{g}\right\} ,
\end{equation*}%
and let $T_{G}^{\ast }Z$ be defined similarly. Let $K_{cpt,G}\left(
T_{G}^{\ast }M\right) $ denote the $G$-equivariant, compactly supported
K-group of $T_{G}^{\ast }M$, which is isomorphic to the group of stable $G$%
-equivariant homotopy classes of transversally elliptic first-order symbols
under direct sum. Likewise, $K_{cpt,H}\left( T^{\ast }Y\right) $ is
isomorphic to the group of the stable $H$-equivariant homotopy classes of
first order elliptic symbols over $Y$.

We define a multiplication%
\begin{equation*}
K_{cpt,G}\left( T_{G}^{\ast }M\right) \otimes K_{cpt,H}\left( T^{\ast
}Y\right) \rightarrow K_{cpt,G}\left( T_{G}^{\ast }Z\right)
\end{equation*}%
as follows. Let $u$ be a transversally elliptic, $G$-equivariant symbol over 
$M$ taking values in $\mathrm{\mathrm{Hom}}\left( E^{+},E^{-}\right) $, and
let $v$ be a $H$-equivariant elliptic symbol over $Y$ taking values in $%
\mathrm{\mathrm{Hom}}\left( F^{+},F^{-}\right) $. First, we lift the symbol $%
u$ to the $H\times G$-equivariant symbol $\widehat{u}$ on $P$. Let $\widehat{%
u}\ast v$ be the standard K-theory multiplication (similar to \cite[Lemma 3.4%
]{A}) 
\begin{equation*}
K_{cpt,H\times G}\left( T_{G}^{\ast }P\right) \otimes K_{cpt,H}\left(
T^{\ast }Y\right) \rightarrow K_{cpt,H\times G}\left( T_{H\times G}^{\ast
}\left( P\times Y\right) \right) .
\end{equation*}%
An element $\left( h,g\right) \in H\times G$ acts on $\left( p,y\right) \in
P\times Y$ by%
\begin{equation}
\left( h,g\right) \left( p,y\right) =\left( phg,h^{-1}y\right) =\left(
pgh,h^{-1}y\right) .  \label{HxGaction}
\end{equation}%
Since the action of $H\times \left\{ e\right\} $ is free, we have%
\begin{equation*}
K_{cpt,H\times G}\left( T_{H\times G}^{\ast }\left( P\times Y\right) \right)
\cong K_{cpt,G}\left( T_{G}^{\ast }\left( P\times _{H}Y\right) \right)
=K_{cpt,G}\left( T_{G}^{\ast }Z\right) .
\end{equation*}%
Finally we define 
\begin{equation*}
u\cdot v=\widetilde{\widehat{u}\ast v}
\end{equation*}%
to be the image of $\widehat{u}\ast v$ in $K_{cpt,G}\left( T_{G}^{\ast
}Z\right) $ under the isomorphism above.

Given any finite-dimensional unitary virtual $H$-representation $\tau $ on $%
V $, we may form the associated $G$-virtual bundle $\widetilde{V_{\tau }}%
=P\times _{\tau }V$ over $M$, defining a class in $K_{G}\left( M\right) $.
The tensor product makes $K_{cpt,G}\left( T_{G}^{\ast }M\right) $ naturally
into a $K_{G}\left( M\right) $-module; for each $\left[ u\right] \in
K_{cpt,G}\left( T_{G}^{\ast }M\right) $, the symbol $u\otimes \tau
:=u\otimes \mathbf{1}_{\widetilde{V_{\tau }}}$ defines an element of $%
K_{cpt,G}\left( T_{G}^{\ast }M\right) $.

We let $\mathrm{ind}^{H}\left( \cdot \right) $ denote the virtual
representation-valued index as explained in \cite{A}; note that the result
is a finite-dimensional virtual representation if the input is a symbol of
an elliptic operator.

\begin{theorem}
\label{EquivariantMultiplicativeTheorem}Let $Z=P\times _{H}Y$ as above, with 
$P$ a $H$-bundle over $M$. Let $u$ be a transversally elliptic, $G$%
-equivariant symbol over $M$ taking values in $\mathrm{\mathrm{Hom}}\left(
E^{+},E^{-}\right) $, and let $v$ be a $H$-equivariant elliptic symbol over $%
Y$ taking values in $\mathrm{\mathrm{Hom}}\left( F^{+},F^{-}\right) $, so
that $u$ and $v$ define classes $\left[ u\right] $ and $\left[ v\right] $ in 
$K_{cpt,G}\left( T_{G}^{\ast }M\right) $ and $K_{cpt,H}\left( T^{\ast
}Y\right) $, respectively. Then $u\cdot v$ defines an element of $%
K_{cpt,G}\left( T_{G}^{\ast }Z\right) $, and 
\begin{equation*}
\mathrm{ind}^{G}\left( u\cdot v\right) =\mathrm{ind}^{G}\left( u\otimes 
\mathrm{ind}^{H}\left( v\right) \right) .
\end{equation*}

\begin{proof}
We adopt the argument in \cite[13.6]{L-M} to our situation. Let $L$ be a
transversally elliptic, $G$-equivariant first order operator representing $u$%
, and let $Q$ be an elliptic, $H$-equivariant first-order operator
representing $v$. Let $\widehat{u}$ be the lift of $u$ to a $H\times G$%
-transversely elliptic symbol over $P$, and let $\widehat{L}$ be a
transversally elliptic, $H\times G$-equivariant first order operator
representing $\widehat{u}$. Next, consider operator product $D=\widehat{L}%
\ast Q$ over $P\times Y$, which represents $\widehat{u}\ast v$. This
operator is $H\times G$ equivariant with respect to the action (\ref%
{HxGaction}). Then%
\begin{equation*}
\ker \left( D^{\ast }D\right) =\left[ \ker \left( \widehat{L}\otimes \mathbf{%
1}\right) \cap \ker \left( \mathbf{1}\otimes Q\right) \right] \oplus \left[
\ker \left( \widehat{L}^{\ast }\otimes \mathbf{1}\right) \cap \ker \left( 
\mathbf{1}\otimes Q^{\ast }\right) \right]
\end{equation*}%
and%
\begin{equation*}
\ker \left( DD^{\ast }\right) =\left[ \ker \left( \widehat{L}^{\ast }\otimes 
\mathbf{1}\right) \cap \ker \left( \mathbf{1}\otimes Q\right) \right] \oplus %
\left[ \ker \left( \widehat{L}\otimes \mathbf{1}\right) \cap \ker \left( 
\mathbf{1}\otimes Q^{\ast }\right) \right]
\end{equation*}%
Let $\widetilde{D}$ and $\widetilde{D}^{\ast }$ be the restrictions of the
operators $D$ and $D^{\ast }$ to sections that are pullbacks of sections
over the base $Z=P\times _{H}Y$, i.e. those that are $H$-invariant. Let $%
\tau ^{+}$ denote the $H$-representation $\ker $ $Q$, and let $\tau ^{-}$ be
the representation $\ker $ $Q^{\ast }$. By the definition of the $H$-action
in (\ref{HxGaction}), the decomposition yields the associated kernels%
\begin{eqnarray*}
\ker \left( \widetilde{D}^{\ast }\widetilde{D}\right) &=&\ker \left( 
\widehat{L}\otimes \tau ^{+}\right) \oplus \ker \left( \widehat{L}^{\ast
}\otimes \tau ^{-}\right) , \\
\ker \left( \widetilde{D}\widetilde{D}^{\ast }\right) &=&\ker \left( 
\widehat{L}^{\ast }\otimes \tau ^{+}\right) \oplus \ker \left( \widehat{L}%
\otimes \tau ^{-}\right) .
\end{eqnarray*}%
We next decompose the above as $G$-representations, and we obtain%
\begin{equation*}
\mathrm{ind}^{G}\left( \widetilde{D}\right) =\mathrm{ind}^{G}\left( \widehat{%
L}\otimes \mathrm{ind}^{H}\left( Q\right) \right) .
\end{equation*}%
The result follows, since $u\cdot v=\widetilde{\widehat{u}\ast v}$ is stably
homotopic to the principal symbol of $\widetilde{D}$.
\end{proof}
\end{theorem}

\subsection{Index of Lifted Dirac operators}

Suppose that $D:\Gamma \left( M,E^{+}\right) \rightarrow \Gamma \left(
M,E^{-}\right) $ is any transversally elliptic, $G$-equivariant operator
with transversally elliptic symbol $u\in \Gamma \left( M,\mathrm{\mathrm{Hom}%
}\left( E^{+},E^{-}\right) \right) $, so that $\left[ u\right] \in
K_{cpt,G}\left( T_{G}^{\ast }M\right) $. If $\mathbf{1}$ denotes the trivial 
$O\left( n\right) $ representation over the identity, then let $\ v$ be any
element of the class $i!\left( \mathbf{1}\right) \in K_{cpt,O\left( n\right)
}\left( T^{\ast }O\left( n\right) \right) $ induced from the inclusion of
the identity in $O\left( n\right) $ via an extension of the Thom isomorphism
(see \cite{ASi1}). Observe that the equivariant index $\mathrm{ind}^{O\left(
n\right) }\left( v\right) $ of the elliptic symbol $v$ is equal to one copy
of the trivial representation (see axioms of the equivariant index in \cite[%
13.6]{ASi1}). By Theorem \ref{EquivariantMultiplicativeTheorem}, $G$%
-equivariant transversally elliptic symbol $u\cdot v$ defines an element of $%
K_{cpt,G}\left( T_{G}^{\ast }F_{O}\right) $ such that%
\begin{eqnarray*}
\mathrm{ind}^{G}\left( u\cdot v\right) &=&\mathrm{ind}^{G}\left( u\otimes 
\mathrm{ind}^{O\left( n\right) }\left( v\right) \right) \\
&=&\mathrm{ind}^{G}\left( u\otimes \mathbf{1}\right) \\
&=&\mathrm{ind}^{G}\left( u\right) .
\end{eqnarray*}

Suppose further that $F_{O}$ is $G$-transversally spin$^{c}$. Then the class 
$\left[ u\cdot v\right] \in K_{cpt,G}\left( T_{G}^{\ast }F_{O}\right) $ may
be represented by the symbol of a transversally-elliptic, $G$-equivariant
operator $D_{N\mathcal{F}}$ of Dirac type.

Thus, the operator $D_{M}^{\mathbf{1}}=i_{\mathbf{1}}^{-1}\circ D_{N\mathcal{%
F}}\circ i_{\mathbf{1}}$ satisfies%
\begin{eqnarray*}
\mathrm{ind}^{G}\left( D_{M}^{\mathbf{1}}\right) &=&\mathrm{ind}^{G}\left(
D_{N\mathcal{F}}\right) =\mathrm{ind}^{G}\left( u\cdot v\right) \\
&=&\mathrm{ind}^{G}\left( u\right) =\mathrm{ind}^{G}\left( D\right) .
\end{eqnarray*}

The result below follows.

\begin{theorem}
Suppose that $F_{O}$ is $G$-transversally spin$^{c}$. Then for every
transversally elliptic symbol class $\left[ u\right] \in K_{cpt,G}\left(
T_{G}^{\ast }M\right) $, there exists an operator of type $D_{M}^{\mathbf{1}%
} $ such that $\mathrm{ind}^{G}\left( u\right) =\mathrm{ind}^{G}\left(
D_{M}^{\mathbf{1}}\right) $.\label{indexClassGivenByTransvDiracThm}
\end{theorem}

\section{Example\label{ExampleSection}}

\subsection{A transversal Dirac operator on the sphere}

Let $G=S^{1}$ act on $S^{2}\subset \mathbb{R}^{3}$ by rotations about the $z$%
-axis. Let $p:F_{O}\rightarrow S^{2}$ be the oriented orthonormal frame
bundle. We will identify $F_{O}$ with $SO\left( 3\right) $ by letting the
first row denote the point on $S^{2}$ and the last two rows denote the
framing of the tangent space. We choose the metric on $F_{O}$ to be 
\begin{equation*}
\left\langle A,B\right\rangle =\mathrm{tr}\left( A^{t}B\right) .
\end{equation*}%
The action of $S^{1}$ lifted to $F_{O}$ is given by multiplication on the
right: 
\begin{equation*}
R_{t}\left( A\right) =A\left( 
\begin{array}{lll}
\cos t & -\sin t & 0 \\ 
\sin t & \cos t & 0 \\ 
0 & 0 & 1%
\end{array}%
\right) .
\end{equation*}%
Tangent vectors to $F_{O}$ are elements of the Lie algebra 
\begin{equation*}
\mathfrak{o}\left( 3\right) =\left\{ \left. \left( 
\begin{array}{lll}
0 & a & b \\ 
-a & 0 & c \\ 
-b & -c & 0%
\end{array}%
\right) \,\right\vert \,a,b,c\in \mathbb{R}\right\} ,
\end{equation*}%
and the tangent space to the $S^{1}$ action is the span of the
left-invariant vector field $T$ induced by $\left( 
\begin{array}{lll}
0 & 1 & 0 \\ 
-1 & 0 & 0 \\ 
0 & 0 & 0%
\end{array}%
\right) $ at the identity. Thus, the normal bundle of the corresponding
foliation on $F_{O}$ is trivial. It is the subbundle $N_{S^{1}}$ of $TF_{O}$
that is given at $A\in F_{O}$ by 
\begin{equation*}
\left. N_{S^{1}}\right\vert _{A}=\left\{ \left. A\left( 
\begin{array}{lll}
0 & 0 & b \\ 
0 & 0 & c \\ 
-b & -c & 0%
\end{array}%
\right) \,\right\vert \,b,c\in \mathbb{R}\right\} .
\end{equation*}%
The vectors $V_{1}=A\left( 
\begin{array}{lll}
0 & 0 & 1 \\ 
0 & 0 & 0 \\ 
-1 & 0 & 0%
\end{array}%
\right) $, $V_{2}=A\left( 
\begin{array}{lll}
0 & 0 & 0 \\ 
0 & 0 & 1 \\ 
0 & -1 & 0%
\end{array}%
\right) $ and the orbit direction $T=A\left( 
\begin{array}{lll}
0 & 1 & 0 \\ 
-1 & 0 & 0 \\ 
0 & 0 & 0%
\end{array}%
\right) $ are mutually orthogonal. Let $E=F_{O}\times \mathbb{C}%
^{2}\rightarrow F_{O}$ be the trivial bundle. The action of $\mathbb{C}%
\mathrm{l}\left( N_{S^{1}}\right) \cong \mathbb{C}\mathrm{l}\left( \mathbb{R}%
^{2}\right) $ on fibers of $E$ is defined by $c\left( V_{1}\right) =\left( 
\begin{array}{ll}
0 & -1 \\ 
1 & 0%
\end{array}%
\right) ,$ $c\left( V_{2}\right) =\left( 
\begin{array}{ll}
0 & i \\ 
i & 0%
\end{array}%
\right) $. We identify the vectors $\left( 
\begin{array}{c}
1 \\ 
0%
\end{array}%
\right) $ and $\left( 
\begin{array}{c}
0 \\ 
1%
\end{array}%
\right) \in \mathbb{C}^{2}$ with the left-invariant fields $V_{1}$ and $%
V_{2} $ in $\Gamma \left( N_{S^{1}}\right) $. We assume that the $S^{1}$%
-action on $E$ is trivial. As in Example \ref{RiemFoliationTrDiracExample},
the transversal Dirac operator is 
\begin{equation*}
D_{N\mathcal{F}}=A_{N\mathcal{F}}=\sum_{j=1}^{2}c\left( V_{j}\right) \nabla
_{V_{j}},
\end{equation*}%
where $\nabla _{V_{j}}$ is the directional derivative in the direction $%
V_{j} $. Since the length of each orbit of the $S^{1}$ action is constant,
the mean curvature vector is zero.

The bundle $F_{O}\rightarrow S^{2}$ is an $SO\left( 2\right) $ principal
bundle and comes equipped with an action of $SO\left( 2\right) $ on the
frames over a point. The left action of%
\begin{equation*}
\left( 
\begin{array}{cc}
\cos \alpha & -\sin \alpha \\ 
\sin \alpha & \cos \alpha%
\end{array}%
\right) \in SO\left( 2\right)
\end{equation*}%
on a frame $A\in SO\left( 3\right) $ is given by%
\begin{equation*}
L_{\alpha }\left( A\right) =\left( 
\begin{array}{ccc}
1 & 0 & 0 \\ 
0 & \cos \alpha & \sin \alpha \\ 
0 & -\sin \alpha & \cos \alpha%
\end{array}%
\right) A.
\end{equation*}%
Again, we extend this action trivially to the $\mathbb{C}^{2}$ bundle. Note
that $D_{N\mathcal{F}}$ is equivariant with respect to both the above left $%
SO\left( 2\right) $ action and the right $S^{1}$ action.

We choose the standard spherical coordinates $x\left( \theta ,\phi \right)
\in S^{2}$. Let $P_{\theta ,\phi }$ denote parallel transport in the tangent
bundle from the north pole along the minimal geodesic connected to $x\left(
\theta ,\phi \right) $. Then%
\begin{equation*}
P_{\theta ,\phi }v=\left( 
\begin{array}{ccc}
\cos \theta & -\sin \theta & 0 \\ 
\sin \theta & \cos \theta & 0 \\ 
0 & 0 & 1%
\end{array}%
\right) \left( 
\begin{array}{ccc}
\cos \phi & 0 & \sin \phi \\ 
0 & 1 & 0 \\ 
-\sin \phi & 0 & \cos \phi%
\end{array}%
\right) \left( 
\begin{array}{ccc}
\cos \theta & \sin \theta & 0 \\ 
-\sin \theta & \cos \theta & 0 \\ 
0 & 0 & 1%
\end{array}%
\right) v.
\end{equation*}%
We parallel transport the standard frame $\left( e_{1},e_{2}\right) $ at the
north pole to get $X_{\theta ,\phi }=P_{\theta ,\phi }e_{1}$, $Y_{\theta
,\phi }=P_{\theta ,\phi }e_{2}$, and the we rotate by $\alpha $ to get all
possible frames. The result is a coordinate chart $U^{1}:\left[ 0,2\pi %
\right] \times \left[ 0,\frac{\pi }{2}\right] \times \left[ 0,2\pi \right]
\rightarrow SO\left( 3\right) $ defined by

\begin{equation*}
U^{1}\left( \theta ,\phi ,\alpha \right) =L_{\alpha }\left( 
\begin{array}{c}
x\left( \theta ,\phi \right) \\ 
X_{\theta ,\phi } \\ 
Y_{\theta ,\phi }%
\end{array}%
\right) .
\end{equation*}

A section $u$ is defined to be of irreducible representation type $\sigma
_{n}:SO\left( 2\right) \rightarrow \mathbb{C}$ if it satisfies $\left(
\left( L_{\beta }\right) _{\ast }u\right) =e^{in\beta }u$. Since the action
of $SO\left( 2\right) $ on the fibers of $E$ is trivial, we have 
\begin{eqnarray}
\left( \left( L_{\beta }\right) _{\ast }u\right) \left( \theta ,\phi ,\alpha
\right) &=&\left( u\circ \left( L_{\beta }\right) ^{-1}\right) \left( \theta
,\phi ,\alpha \right)  \notag \\
&=&u\left( \theta ,\phi ,\alpha -\beta \right) =e^{in\beta }u\left( \theta
,\phi ,\alpha \right) .  \label{leftMultiplicationFormula}
\end{eqnarray}%
Thus, it suffices to calculate $u\left( \theta ,\phi ,0\right) =u$ at $%
U^{1}\left( \theta ,\phi ,0\right) $.

The lower hemisphere coordinates of the point and vectors would have the
opposite third coordinate, and the sign of $X_{\theta ,\phi }$ is reversed
in addition to ensure that the frame is oriented. Note that the $\phi $ in
the lower hemisphere is $\left( \pi -\phi \right) $ in the upper hemisphere.
Thus the second chart is $U^{2}:\left[ 0,2\pi \right] \times \left[ 0,\frac{%
\pi }{2}\right] \times \left[ 0,2\pi \right] $%
\begin{equation*}
U^{2}\left( \theta ,\phi ,\alpha \right) =U_{\alpha }^{2}\left( 
\begin{array}{c}
x\left( \theta ,\phi \right) \\ 
-X_{\theta ,\phi } \\ 
Y_{\theta ,\phi }%
\end{array}%
\right) \mathrm{diag}\left( 1,1,-1\right)
\end{equation*}%
One can check that 
\begin{equation*}
U^{1}\left( \theta ,\frac{\pi }{2},\alpha \right) =U^{2}\left( \theta ,\frac{%
\pi }{2},\alpha -2\theta \right) =L_{-2\theta }U^{2}\left( \theta ,\frac{\pi 
}{2},\alpha \right) .
\end{equation*}%
Thus, the clutching function for the frame bundle is multiplication on the
left by $e^{2\theta i}$.

Next, suppose that $u$ is a section such that $\left( \left( L_{\beta
}\right) _{\ast }u\right) \left( M\right) =u\left( L_{-\beta }M\right)
=e^{in\beta }u\left( M\right) $. This means in fact that $u\left( \theta
,\phi ,\alpha -\beta \right) =e^{in\beta }u\left( \theta ,\phi ,\alpha
\right) $ in both charts. Thus we may trivialize the bundle by restricting
to $\left( \theta ,\phi ,0\right) $ in each chart. We observe

\begin{equation*}
u^{1}\left( \theta ,\frac{\pi }{2},0\right) =u^{2}\left( L_{-2\theta }\left(
\theta ,\frac{\pi }{2},0\right) \right) =e^{i2n\theta }u^{2}\left( \theta ,%
\frac{\pi }{2},0\right) ,
\end{equation*}%
and thus the clutching function for $\mathbb{E}^{\sigma _{n}}$ is $%
e^{2n\theta i}$.

One may express the vector fields $V_{1}$, $V_{2}$ in terms of the
coordinate vector fields $\partial _{\alpha }$,$\partial _{\theta }$, $%
\partial _{\phi }$. In the upper hemisphere,%
\begin{eqnarray*}
V_{1}^{1} &=&\frac{\sin \theta \left( \cos \phi -1\right) }{\sin \phi }%
\partial _{\alpha }+\frac{\sin \theta \cos \phi }{\sin \phi }\partial
_{\theta }-\cos \theta ~\partial _{\phi } \\
V_{2}^{1} &=&\frac{\cos \theta \left( 1-\cos \phi \right) }{\sin \phi }%
\partial _{\alpha }-\frac{\cos \theta \cos \phi }{\sin \phi }\partial
_{\theta }-\sin \theta ~\partial _{\phi }
\end{eqnarray*}

Now we wish to consider the operator%
\begin{equation*}
D_{S^{2}}^{\sigma _{n}}=i_{\sigma _{n}}^{-1}\circ D_{N\mathcal{F}}\circ
i_{\sigma _{n}}:\Gamma \left( S^{2},\mathcal{E}^{\sigma _{n}}\right)
\rightarrow \Gamma \left( S^{2},\mathcal{E}^{\sigma _{n}}\right) ,
\end{equation*}%
where $i_{\sigma _{n}}:\Gamma \left( S^{2},\mathcal{E}^{\sigma _{n}}\right)
\rightarrow \Gamma \left( F_{O},E\right) ^{\sigma _{n}}$. We have%
\begin{eqnarray*}
\sqrt{2}D_{N\mathcal{F}}^{1}\left( 
\begin{array}{c}
u_{1} \\ 
u_{2}%
\end{array}%
\right) &=&\left( 
\begin{array}{ll}
0 & -1 \\ 
1 & 0%
\end{array}%
\right) \nabla _{V_{1}^{1}}\left( 
\begin{array}{c}
u_{1} \\ 
u_{2}%
\end{array}%
\right) +\left( 
\begin{array}{ll}
0 & i \\ 
i & 0%
\end{array}%
\right) \nabla _{V_{2}^{1}}\left( 
\begin{array}{c}
u_{1} \\ 
u_{2}%
\end{array}%
\right) \\
&=&\left( 
\begin{array}{c}
-\overline{\left( V_{1}^{1}+iV_{2}^{1}\right) }u_{2} \\ 
\left( V_{1}^{1}+iV_{2}^{1}\right) u_{1}%
\end{array}%
\right) .
\end{eqnarray*}%
There is a similar formula for $\sqrt{2}D_{N\mathcal{F}}^{2}$ in the lower
hemisphere chart. Observe that%
\begin{equation*}
\left( V_{1}^{1}+iV_{2}^{1}\right) =-ie^{i\theta }\left( \cot \phi -\csc
\phi \right) \partial _{\alpha }+-ie^{i\theta }\cot \phi \partial _{\theta
}+-e^{i\theta }\partial _{\phi }
\end{equation*}

We easily check that right multiplication by $\beta \in S^{1}=\mathbb{R~}%
\mathrm{mod}\,2\pi $ on $U^{1}\left( \theta ,\phi ,\alpha \right) $ satisfies

\begin{equation*}
R_{\beta }U^{1}\left( \theta ,\phi ,\alpha \right) =U^{1}\left( \theta
+\beta ,\phi ,\alpha +\beta \right) .
\end{equation*}

If $\psi ^{1}\left( \theta ,\phi ,0\right) $ is a section of $F_{O}\times 
\mathbb{C}^{2}\rightarrow F_{O}$ of type $\sigma _{n}$ (with respect to the
fiberwise action of $SO\left( 2\right) $) over the upper hemisphere, then%
\begin{eqnarray*}
\left( \left( R_{\beta }\right) _{\ast }\psi ^{1}\right) \left( \theta ,\phi
,0\right) &=&\psi ^{1}\circ R_{-\beta }\left( \theta ,\phi ,0\right) \\
&=&\psi ^{1}\left( \theta -\beta ,\phi ,-\beta \right) =e^{in\beta }\psi
^{1}\left( \theta -\beta ,\phi ,0\right) ,
\end{eqnarray*}%
using the upper hemisphere trivialization $U^{1}\left( \theta ,\phi ,\alpha
\right) $ and equation (\ref{leftMultiplicationFormula}).

\vspace{0in}If we assume that $\psi :F_{O}\rightarrow \mathcal{E}^{\sigma
_{n}}$ is a section of type $\rho _{m}$ with respect to the lifted $S^{1}$
action, then $\left( \left( R_{\beta }\right) _{\ast }\psi ^{1}\right)
=e^{im\beta }\psi ^{1}$. Thus, 
\begin{equation*}
\left( \left( R_{\beta }\right) _{\ast }\psi ^{1}\right) \left( \theta ,\phi
,0\right) =e^{im\beta }\psi ^{1}\left( \theta ,\phi ,0\right) =e^{in\beta
}\psi ^{1}\left( \theta -\beta ,\phi ,0\right) ,
\end{equation*}%
which implies%
\begin{equation*}
\psi ^{1}\left( \theta ,\phi ,0\right) =e^{i\left( n-m\right) \theta }\psi
^{1}\left( 0,\phi ,0\right) .
\end{equation*}%
The analogous calculation in the lower hemisphere chart yields%
\begin{equation*}
\psi ^{2}\left( \theta ,\phi ,0\right) =e^{i\left( -n-m\right) \theta }\psi
^{2}\left( 0,\phi ,0\right) .
\end{equation*}%
\vspace{0in}

\subsection{Calculation of $\ker $ $D_{S^{2}}^{\protect\sigma _{n}}$}

Since $D_{S^{2}}^{\sigma _{n}}=i_{\sigma _{n}}^{-1}\circ D_{N\mathcal{F}%
}\circ i_{\sigma _{n}}$, we seek solutions to the equation%
\begin{equation*}
\sqrt{2}D_{N\mathcal{F}}^{1}\left( 
\begin{array}{c}
\psi _{1} \\ 
\psi _{2}%
\end{array}%
\right) =\left( 
\begin{array}{c}
-\overline{\left( V_{1}^{1}+iV_{2}^{1}\right) }\psi _{2} \\ 
\left( V_{1}^{1}+iV_{2}^{1}\right) \psi _{1}%
\end{array}%
\right) =\left( 
\begin{array}{c}
0 \\ 
0%
\end{array}%
\right)
\end{equation*}%
in the upper hemisphere chart. From the equations $\left(
V_{1}^{1}+iV_{2}^{1}\right) \psi _{1}=0$, $\partial _{\alpha }\psi
_{1}=-ni\psi _{1}$ (since $\Psi \in \Gamma \left( S^{2},\mathcal{E}^{\sigma
_{n}}\right) $ ), and $\partial _{\theta }\psi _{1}=i\left( n-m\right) \psi
_{1}$ (since $\Psi \in \Gamma \left( S^{2},\mathcal{E}^{\sigma _{n}}\right)
^{\rho _{m}}$), we have%
\begin{eqnarray*}
0 &=&\left( V_{1}^{1}+iV_{2}^{1}\right) \psi _{1}=\left( -ie^{i\theta
}\left( \cot \phi -\csc \phi \right) \partial _{\alpha }+-ie^{i\theta }\cot
\phi \partial _{\theta }+-e^{i\theta }\partial _{\phi }\right) \psi _{1} \\
&=&\left( -me^{i\theta }\left( \cot \phi \right) +ne^{i\theta }\left( \csc
\phi \right) -e^{i\theta }\partial _{\phi }\right) \psi _{1}.
\end{eqnarray*}%
Solving this equation, we obtain%
\begin{equation*}
\psi _{1}\left( 0,\phi \right) =C_{2}\frac{\left( \sin \phi \right) ^{n-m}}{%
\left( \cos \phi +1\right) ^{n}}.
\end{equation*}%
This implies that $\psi _{1}$ is%
\begin{eqnarray*}
\psi _{1}\left( \theta ,\phi \right) &=&C_{2}\frac{\left( \sin \phi \right)
^{n-m}}{\left( \cos \phi +1\right) ^{n}}e^{i\left( n-m\right) \theta },\text{
or} \\
\psi _{1}\left( z\right) &=&C_{2}\frac{z^{n-m}}{\left( \sqrt{1-\left\vert
z\right\vert ^{2}}+1\right) ^{n}}\text{ }
\end{eqnarray*}%
in the complex coordinates of the projection of the upper hemisphere to the $%
xy$ plane. Thus, $\psi _{1}$ is smooth in the upper hemisphere only if $%
n\geq m$. Similarly, $-\overline{\left( V_{1}^{1}+iV_{2}^{1}\right) }\psi
_{2}=0$, $\partial _{\alpha }\psi _{2}=-ni\psi _{2}$, and $\partial _{\theta
}\psi _{2}=i\left( n-m\right) \psi _{2}$ implies%
\begin{eqnarray*}
\psi _{2}\left( \theta ,\phi \right) &=&C\left( \cos \phi +1\right)
^{n}\left( \sin \phi \right) ^{m-n}e^{i\left( n-m\right) \theta },\text{ or}
\\
\psi _{2}\left( z\right) &=&C_{2}\left( \sqrt{1-\left\vert z\right\vert ^{2}}%
+1\right) ^{n}\overline{z}^{m-n}.
\end{eqnarray*}%
Hence, $\psi _{2}$ is smooth in the upper hemisphere only if $m\geq n$.

\vspace{0in}We need to see if the solutions $\psi _{1}$ and $\psi _{2}$
extend to solutions over the entire sphere. In the lower hemisphere, we have
the equation%
\begin{equation*}
\sqrt{2}D_{N\mathcal{F}}^{2}\Psi =\sqrt{2}D_{N\mathcal{F}}^{2}\left( 
\begin{array}{c}
\psi _{1} \\ 
\psi _{2}%
\end{array}%
\right) =\left( 
\begin{array}{c}
\left( -V_{1}^{2}+iV_{2}^{2}\right) \psi _{2} \\ 
\left( V_{1}^{2}+iV_{2}^{2}\right) \psi _{1}%
\end{array}%
\right) =\left( 
\begin{array}{c}
0 \\ 
0%
\end{array}%
\right) .
\end{equation*}%
Similar computations show that%
\begin{eqnarray*}
\psi _{1}\left( z\right) &=&C_{2}\left( \sqrt{1-\left\vert z\right\vert ^{2}}%
+1\right) ^{n}z^{-n-m}, \\
\psi _{2}\left( z\right) &=&C_{2}\left( \sqrt{1-\left\vert z\right\vert ^{2}}%
+1\right) ^{-n}\overline{z}^{m+n}
\end{eqnarray*}%
in the complex coordinates of the projection of the lower hemisphere to the $%
xy$ plane. Thus, $\psi _{1}$ is smooth in the lower hemisphere only if $%
n+m\leq 0$, and $\psi _{2}$ is smooth in the lower hemisphere only if $%
m+n\geq 0$.

In summary, we seek solutions to $\sqrt{2}D_{N\mathcal{F}}^{2}\Psi =\sqrt{2}%
D_{N\mathcal{F}}^{1}\Psi =0$ restricted to sections of $\mathcal{E}^{\sigma
_{n}}$ of type $\rho _{m}$. The clutching function of $\mathcal{E}^{\sigma
_{n}}$ is multiplication by $e^{2n\theta i}$ (i.e. $z^{2n}$ or $\overline{z}%
^{-2n}$), so that $\Psi ^{1}\left( \theta ,\frac{\pi }{2},0\right) =\Psi
^{2}\left( A_{2\theta }\left( \theta ,\frac{\pi }{2},0\right) \right)
=e^{i2n\theta }\Psi ^{2}\left( \theta ,\frac{\pi }{2},0\right) $.

The function $\psi _{1}$ is continuous if and only if $\psi _{1}^{1}\left(
\theta ,\frac{\pi }{2},0\right) =e^{i2n\theta }\psi _{1}^{2}\left( \theta ,%
\frac{\pi }{2},0\right) $ and $m\leq -\left\vert n\right\vert $. Similarly, $%
\psi _{2}$ is continuous if and only if $\psi _{2}^{1}\left( \theta ,\frac{%
\pi }{2},0\right) =e^{i2n\theta }\psi _{2}^{2}\left( \theta ,\frac{\pi }{2}%
,0\right) $ and $m\geq \left\vert n\right\vert $.

\vspace{0in}From the equations in the previous section, the index of $%
D_{S^{2}}^{\sigma _{n}}$ restricted to sections of type $\rho _{m}$ is 
\begin{equation}
\mathrm{ind}^{\rho _{m}}\left( D_{S^{2}}^{\sigma _{n}}\right) =\left\{ 
\begin{array}{cc}
-1 & \text{if }m>\left\vert n\right\vert \text{ or }m=\left\vert
n\right\vert ~\text{and }n\neq 0 \\ 
0 & \text{if }-\left\vert n\right\vert <m<\left\vert n\right\vert \text{ or }%
m=n=0 \\ 
1 & \text{if }m<-\left\vert n\right\vert \text{ or }m=-\left\vert
n\right\vert ~\text{and }n\neq 0%
\end{array}%
\right. .  \label{indexFormula}
\end{equation}%
Note that the kernel of $D_{S^{2}}^{\sigma _{n}}$ is infinite-dimensional.
The operator $D_{S^{2}}^{\sigma _{n}}$ fails to be elliptic precisely at the
points where $\cot \phi =0$; that is, at the equator.

\subsection{The operator on $F_{O}\diagup G$}

We now construct the operator $D_{F_{O}\diagup G}^{\rho _{m}}:\Gamma \left(
F_{O}\diagup G,\mathcal{T}^{\rho _{m}}\right) \rightarrow \Gamma \left(
F_{O}\diagup G,\mathcal{T}^{\rho _{m}}\right) $. First, observe that $%
F_{O}\diagup G=SO\left( 3\right) \diagup S^{1}$ is again the sphere $S^{2}$.
The orbits of the action of $G$ on $F_{O}$ are of the form $\left\{
R_{t}M:t\in \left[ 0,2\pi \right] \right\} $, with $M\in F_{O}$. The map $%
R_{t}$ rotates the first and second columns of the matrix, so that the map $%
F_{O}\rightarrow S^{2}$ is the map to the third column. Thus, the projection
of the vector $V_{1}$ to $TS^{2}$ is

\begin{equation*}
V_{1}=M\left( 
\begin{array}{lll}
0 & 0 & 1 \\ 
0 & 0 & 0 \\ 
-1 & 0 & 0%
\end{array}%
\right) \mapsto \text{first column of }M
\end{equation*}

Similarly,%
\begin{equation*}
V_{2}=M\left( 
\begin{array}{lll}
0 & 0 & 0 \\ 
0 & 0 & 1 \\ 
0 & -1 & 0%
\end{array}%
\right) \mapsto \text{second column of }M
\end{equation*}

A section of type $\rho _{m}$ of $F_{O}\times \mathbb{C}^{2}\rightarrow
F_{O} $ is one for which the partial derivative in direction 
\begin{equation*}
T=M\left( 
\begin{array}{lll}
0 & 1 & 0 \\ 
-1 & 0 & 0 \\ 
0 & 0 & 0%
\end{array}%
\right)
\end{equation*}%
is multiplication by $im$. In the $U^{1}$ coordinate chart (corresponding to 
$0\leq \phi \leq \frac{\pi }{2}$, $0\leq \theta \leq 2\pi $, $0\leq \alpha
\leq 2\pi $), the quotient to $F_{O}\diagup G$ goes to $\left( 
\begin{array}{c}
\cos \phi \\ 
-\cos \left( \theta -\alpha \right) \sin \phi \\ 
-\sin \left( \theta -\alpha \right) \sin \phi%
\end{array}%
\right) $, mapping to the entire upper hemisphere $x\geq 0$, with fibers of
the form $\left( \phi ,\theta ,\alpha \right) =\left( \phi _{0},\theta
_{0}+t,\alpha _{0}+t\right) $. Thus, we may fix $\alpha =0$ for the sake of
argument and allow $\theta =\theta -\alpha $ and $\phi $ to vary. The group $%
S^{1}$ acts as before by $U^{1}\left( \theta ,\phi ,\alpha \right) \overset{%
R_{\beta }}{\mapsto }U^{1}\left( \theta +\beta ,\phi ,\alpha +\beta \right) $%
. In particular, if $\psi $ is a section of type $\rho _{m}$, $\left(
R_{\beta }\psi ^{1}\right) \left( \theta ,\phi ,\alpha \right) =\psi
^{1}\circ m_{-\beta }\left( \theta ,\phi ,\alpha \right) =\psi ^{1}\left(
\theta -\beta ,\phi ,\alpha -\beta \right) =e^{im\beta }\psi ^{1}\left(
\theta ,\phi ,\alpha \right) $. Thus, setting $\alpha =0$, we have%
\begin{equation*}
\psi ^{1}\left( \theta -\beta ,\phi ,-\beta \right) =e^{im\beta }\psi \left(
\theta ,\phi ,0\right) .
\end{equation*}%
Thus,%
\begin{eqnarray*}
\frac{\partial }{\partial \beta }\left[ \psi ^{1}\left( \theta -\beta ,\phi
,-\beta \right) \right] &=&im\psi ^{1}\left( \theta -\beta ,\phi ,-\beta
\right) \\
&=&\left( \left( -\frac{\partial }{\partial \theta }-\frac{\partial }{%
\partial \alpha }\right) \psi ^{1}\right) \left( \theta -\beta ,\phi ,-\beta
\right)
\end{eqnarray*}%
So all sections of type $\rho _{m}$ satisfy $\left( -\frac{\partial }{%
\partial \theta }-\frac{\partial }{\partial \alpha }\right) \psi ^{1}=im\psi
^{1},$ or $\partial _{\alpha }=-\partial _{\theta }-im$.

Restricted to this space of sections, we have%
\begin{equation*}
\left( V_{1}^{1}+iV_{2}^{1}\right) =-ie^{i\theta }\csc \phi \partial
_{\theta }-e^{i\theta }\partial _{\phi }-me^{i\theta }\left( \cot \phi -\csc
\phi \right)
\end{equation*}%
The interested reader may check that this operator is elliptic at all points
of the hemisphere $x\geq 0$. A similar statement is true in the other
hemisphere. Thus, $D_{F_{O}\diagup G}^{\rho _{m}}:\Gamma \left( F_{O}\diagup
G,\mathcal{T}^{\rho _{m}}\right) \rightarrow \Gamma \left( F_{O}\diagup G,%
\mathcal{T}^{\rho _{m}}\right) $ is elliptic, as expected.

The kernel of $D_{F_{O}\diagup G}^{\rho _{m}}$ restricted to $\Gamma \left(
F_{O}\diagup G,\mathcal{T}^{\rho _{m}}\right) ^{\sigma _{n}}$ has dimension 
\begin{eqnarray*}
\dim \left( \ker ^{\sigma _{n}}\left( D_{F_{O}\diagup G}^{\rho _{m}}\right)
\right) &=&\dim \left( \ker ^{\rho _{m}}\left( D_{S^{2}}^{\sigma
_{n}}\right) \right) \\
&=&\left\{ 
\begin{array}{ll}
2~ & \text{if }m=n=0 \\ 
1 & \text{if }\left\vert n\right\vert \leq \left\vert m\right\vert ,~m\neq 0
\\ 
0 & \text{if }\left\vert n\right\vert >\left\vert m\right\vert%
\end{array}%
\right. ~,
\end{eqnarray*}%
using the results preceding formula (\ref{indexFormula}). As expected, for a
given bundle $\mathcal{T}^{\rho _{m}}$, only a finite number of the
representation types $\sigma _{n}$ occur, and the kernel of $D_{F_{O}\diagup
G}^{\rho _{m}}$ on the space of all sections is finite-dimensional.


\begin{thebibliography}{99}
\bibitem{A} M. F. Atiyah, \textit{Elliptic operators and compact groups},
Lecture Notes in Math. \textbf{401}, Springer-Verlag, Berlin, 1974.

\bibitem{ASe} M. F. Atiyah and G. B. Segal, \emph{The index of elliptic
operators: II}, Ann. of Math. (2) \textbf{87}(1968), 531--545.

\bibitem{ASi1} M. F. Atiyah and I. M. Singer, \emph{The index of elliptic
operators I}, Ann. of Math. (2) \textbf{87}(1968), 484--530.

\bibitem{B-G-V} N. Berline, E. Getzler, and M. Vergne, \textit{Heat kernels
and Dirac operators}, Grundlehren der Mathematischen Wissenschaften \textbf{%
\ 298}, Springer-Verlag, Berlin, 1992.

\bibitem{Be-V1} N. Berline and M. Vergne, \textit{The Chern character of a
transversally elliptic symbol and the equivariant index}, Invent. Math. 
\textbf{124}(1996), no. 1-3, 11-49.

\bibitem{Be-V2} N. Berline and M. Vergne, \textit{L'indice \'{e}quivariant
des op\'{e}rateurs transversalement elliptiques}, Invent. Math. \textbf{124}
(1996), no. 1-3, 51-101.

\bibitem{Bra3} M. Braverman, \textit{\ Index theorem for equivariant Dirac
operators on noncompact manifolds}, K-Theory \textbf{27} (2002), no. 1,
61--10.

\bibitem{BrKRi} J. Br\"{u}ning, F. W. Kamber, and K. Richardson, \emph{The
eta invariant and equivariant index of transversally elliptic operators},
preprint in preparation.

\bibitem{C} P. R. Chernoff, \textit{Essential self-adjointness of powers of
generators of hyperbolic equations}. J. Functional Analysis \textbf{12}
(1973), 401--414.

\bibitem{DouglasGlK} R. G. Douglas, J. F. Glazebrook, and F. W. Kamber, and
G. L. Yu, \emph{Index formulas for geometric Dirac operators in Riemannian
foliations}, K-Theory \textbf{9}(1995), no. 5, 407--441.

\bibitem{GlK} J. F. Glazebrook and F. W. Kamber, \emph{Transversal Dirac
families in Riemannian foliations}, Comm. Math. Phys. \textbf{140} (1991),
no. 2, 217--240.

\bibitem{Hab} G. Habib, \emph{Eigenvalues of the transversal Dirac operator
on K\"{a}hler foliations}, J. Geom. Phys. 5\textbf{6} (2005), 260--270.

\bibitem{Jung1} S. D. Jung, \emph{The first eigenvalue of the transversal
Dirac operator}, J. Geom. Phys. \textbf{39} (2001), 253--264.

\bibitem{Jung2} S. D. Jung, D.S. Kang, \emph{Eigenvalues of the basic Dirac
operator and transversal twistor operator on Riemannian foliations}, Far
East J. Math. Sci. \textbf{14} (2004), 361--369.

\bibitem{La} C. Lazarov, \emph{Transverse index and periodic orbits}, Geom.
Funct. Anal. \textbf{10} (2000), no. 1, 124--159.

\bibitem{L-M} H. B. Lawson, Jr. and M-L. Michelsohn, \textit{Spin Geometry}
, Princeton Mathematical Series \textbf{38}, Princeton University Press,
Princeton, 1989.

\bibitem{Paradan} P.-E. Paradan, \emph{Spin\_c quantization and the
K-multiplicities of the discrete series}, Ann. Sci. \'{E}cole Norm. Sup. (4) 
\textbf{36} (2003), no. 5, 805--845.
\end{thebibliography}
\end{document}